\documentclass[11pt]{article}
\usepackage{amssymb}
\usepackage{amsmath}
\usepackage[dvips]{graphics}
\usepackage{epsfig}

\begin{document}
\date{}
\author{\textbf{Vassilis G. Papanicolaou}
\\\\
Department of Mathematics
\\
National Technical University of Athens,
\\
Zografou Campus, 157 80, Athens, GREECE
\\\\
{\tt papanico@math.ntua.gr}}
\title{Asymptotics Related to a Binary Search Scheme}
\maketitle
\begin{abstract}
Specimens are collected from $N$ different sources. Each specimen has probability $p$ of being contaminated, independently of the other specimens. We assume group 
testing is applicable, namely one can take small portions from several
specimens, mix them together, and test the mixture for contamination, so that if the test turns positive, then at least one of the samples in the
mixture is contaminated.

In this paper we derive asymptotics, as $N$ gets large, of the expectation and the variance of the number $T(N)$ of tests required in order to find all contaminated 
specimens, under the binary search scheme we introduced in \cite{P} (see, also, arXiv:2007.11910). In \cite{P} the probability $p$ was fixed, whereas in the present
work we consider the case where $p \sim a / N^{-\beta}$, $a, \beta > 0$ with emphasis on the case $\beta = 1$, which turns out to be the most interesting.
\end{abstract}
\textbf{Keywords.} Prevalence (rate); adaptive group testing; binary search scheme; probabilistic testing;
average-case aspect ratio; asymptotics.\\\\
\textbf{2020 AMS Mathematics Subject Classification.} 60C99; 62P10; 92C50.

\section{Introduction}
Consider $N$ containers containing samples from $N$ different sources (e.g., clinical specimens from $N$
different patients, water samples from $N$ different lakes, soil specimens,
food specimens, etc.). For each sample we assume that the probability of being contaminated (by, say, a
virus, a toxic substance, etc.) is $p$ and the probability that it is not contaminated is $q := 1-p$, independently of the other samples.
For instance, in the case of an infectious disease, $p$ is the so-called \emph{prevalence rate}. All $N$ samples must undergo a screening
procedure (say a molecular or antibody test in the case of a viral contamination, a radiation measurement in the case of a radioactive contamination, etc.) in order to 
determine exactly which are the contaminated ones.

One obvious way to identify all contaminated specimens is the \emph{individual testing}, namely to test the specimens one by one individually
(of course, this approach requires $N$ tests). In \cite{P} (see, also, arXiv:2007.11910) we analyzed a group testing (or pool testing)
approach, which we called ``binary search scheme" requiring a random number of tests. In particular, we showed that
if $p$ is small (actually, if $p < 0.224$), then by following
this scheme, the expected number of tests required in order to determine all contaminated samples is strictly less than $N$.

Group testing procedures have broad applications in areas ranging from medicine and engineering to airport security control, and they have attracted
the attention of researchers for over eighty years (see, e.g., Dorfman's seminal paper \cite{Do}, as well as \cite{H}, \cite{S-G}, \cite{U},
etc.). In fact, recently, group testing has received some special attention, partly due to the COVID-19 crisis
(see, e.g., \cite{A1}, \cite{A2}, \cite{A-F-F-D}, \cite{G-G}, \cite{Ma}, \cite{M}, and the references therein).

Our binary search scheme is quite natural and easy to implement. It goes as follows: First we take samples from each of the $N$ containers, mix them together and test the mixture. If the mixture is not contaminated, we know that none of the samples is contaminated. If the mixture is contaminated,
we split the containers in two groups, where the first group consists of the first $\lfloor N/2 \rfloor$ containers
and the second group consists of the remaining $\lceil N/2 \rceil$ containers. Then we take samples from the containers of the first group,
mix them together and test the mixture. If the mixture is not contaminated, then we know that none of the samples of the first group is contaminated.
If the mixture is contaminated, we split the containers of the first group into two subgroups (in the same way we did with the $N$ containers) and
continue the procedure. We also apply the same procedure to the second group (consisting of $\lceil N/2 \rceil$ containers).

Another instance where the above procedure applies is the following. Suppose we have a collection of $N$ identically
looking coins (or cans, or boxes, or pills, etc.), in which a relatively small percentage $p$ of them are counterfeit (defective). All genuine coins
have the same weight, say $w_0$, while the weight of a counterfeit is $< w_0$. If we are interested in finding all counterfeits, we can, again,
apply the binary search scheme. First we weigh all coins simultaneously. If their weight is $N w_0$, we know there are not any
counterfeits. If, however, their weight is $< N w_0$, then there are counterfeits and, as in the case of the
contaminated samples, in order to find them we continue the search by splitting the collection into two subcollections consisting of
$\lfloor N/2 \rfloor$ and $\lceil N/2 \rceil$ coins respectively.

Naturally, the main quantity of interest is the number $T(N) = T(N; q)$ of tests required in order to determine all contaminated samples,
if the above binary search scheme is applied. Following the notation of \cite{P}, in the case where $N = 2^n$ we use $W_n$ to denote the random variable $T(2^n)$.

\medskip

We can characterize the aforementioned procedure as an
\emph{adaptive probabilistic group testing} \cite{A1}.
``Adaptive" refers to the fact that the way we perform a stage of the testing procedure depends on the outcome of the previous
stages, while ``probabilistic" refers to the fact that each specimen has probability $p$ of being contaminated, independently of the other specimens and,
furthermore, that the required number of tests, $T(N)$, is a random variable. Let us be a bit more specific about the latter. Suppose
$\xi_1, \xi_2, \ldots$ is a sequence of independent Bernoulli random variables of a probability space $(\Omega, \mathcal{F}, \mathbb{P})$, with
$\mathbb{P}\{\xi_k = 1\} = p$ and $\mathbb{P}\{\xi_k = 0\} = q$, $k = 1, 2, \ldots$, so that $\{\xi_k = 1\}$ is the event that the $k$-th specimen is contaminated. Then, for each $N$, the quantity $T(N)$ of our binary search scheme is a specific function of $\xi_1, \xi_2, \ldots, \xi_N$, say,
\begin{equation}
T(N) = T_N(\xi_1, \xi_2, \ldots, \xi_N).
\label{I1}
\end{equation}
In particular, is not hard to see that
\begin{equation}
1 = T_N(0, 0, \ldots, 0) \leq T_N(\xi_1, \xi_2, \ldots, \xi_N) \leq T_N(1, 1, \ldots, 1) = 2N - 1,
\label{I2}
\end{equation}
i.e. if none of the contents of the $N$ containers is contaminated,
then $T(N) = 1$, whereas if all $N$ containers contain contaminated samples, then by easy induction on $N$ we can see that $T(N) = 2N-1$
(observe that, for $N \geq 2$, if $\xi_1 = \cdots = \xi_N = 1$ then $T(N) = T(\lfloor N/2 \rfloor) + T(\lceil N/2 \rceil) + 1$, while, of course, $T(1) = 1$).
Thus, $T(N)$ can become bigger than $N$, while the deterministic way of checking of the samples one by one
requires $N$ tests.

Let us also notice that
\begin{equation}
T_N(\xi_1, \xi_2, \ldots, \xi_N) \leq  T_{N+1}(\xi_1, \xi_2, \ldots, \xi_N, \xi_{N+1})
\label{I3}
\end{equation}
(if $\xi_1 = \xi_2 = \cdots = \xi_N = \xi_{N+1} = 0$, then both quantities are equal to $1$; in all other cases we can prove the inequality by induction, by assyming
it is true for all $N < M$ and showing that it also holds for $N = M$).
Evidently, \eqref{I3} implies that
\begin{equation}
T(N) \leq_{st} T(N+1),
\label{I4}
\end{equation}
where $\leq_{st}$ denotes the usual stochastic ordering
(recall that $X \leq_{st} Y$ means that $\mathbb{P}\{X > x\} \leq \mathbb{P}\{Y > x\}$ for all $x \in \mathbb{R}$).
Also, it follows easily from a coupling argument that
\begin{equation}
q_1 > q_2 \ \ \Longrightarrow \ \ T(N; q_1) \leq_{st} T(N; q_2).
\label{I5}
\end{equation}

\medskip

In \cite{P} we derived explicit formulas for the expectation and the variance of $W_n = T(2^n)$ as well as the asymptotic behavior of all the moments of $W_n$,
as $n \to \infty$, assuming that $q$ stays fixed. From the moment asymptotics we concluded that, as $n \to \infty$, the quantity $2^{-n/2} [W_n - 2^n \alpha_1(q)]$, 
where $\alpha_1(q)$ is a specific deterministic quantity, converges in distribution to a zero-mean normal random variable. In the case of an arbitrary $N$ we derived recursive formulas for $\mu(N) := \mathbb{E}[T(N)]$, $g(z; N) := \mathbb{E}[z^{T(N)}]$, $\sigma^2(N) := \mathbb{V}[T(N)]$,
and we showed that the moments of $Y(N) := [T(N) - \mu(N)] / \sigma(N)$ converge to the moments of the standard normal variable $Z$ as $N \to \infty$,
with an immediate consequence that $Y(N)$ converges to $Z$ in distribution. Furthermore, we demonstrated that the limit of
$\mathbb{E}[T(N)] / N$, as $N \to \infty$, does not exist, which is in contrast to the special case $N = 2^n$, where the limit exists.

In the present work we consider the case where $N \to \infty$, while $p \sim a N^{-\beta}$ for some $a, \beta > 0$, thus $p \to 0$. The most interesting case,
at least theoretically, is the case where $p \sim a / N$ (i.e. $\beta = 1$). Notice that in this case the number of contaminated specimens follows (in the limit)
a Poisson distribution with parameter $a$. Our results demonstrate that in the aforementioned regimes our binary search scheme behaves much better than individual testing.

\section{Asymptotics of the expectation and the variance of $T(N)$}
We start with a well-known result \cite{Pe}, which we present here in some detail for the sake of completeness.

\medskip

\textbf{Lemma 1.}
\begin{equation}
1 - \gamma x \leq (1 - x)^{\gamma}
\qquad
\text{whenever}\ \ x \leq 1 \ \ \text{and}\ \ \gamma \geq 1,
\label{A1}
\end{equation}
with equality only in the case where $x = 0$ or $\gamma = 1$.

Also,
\begin{equation}
(1 - x)^{\gamma} \leq 1 - \gamma x + \frac{\gamma(\gamma - 1)}{2} x^2
\qquad
\text{whenever}\ \ x \in [0, 1] \ \ \text{and}\ \ \gamma \geq 2,
\label{A2}
\end{equation}
with equality only in the case where $x = 0$ or $\gamma = 2$.

\smallskip

\textit{Proof}. Let us give the proof of \eqref{A2}, since inequality \eqref{A1} is a version of the Bernoulli inequality and it can be proved in a similar fashion.

If $\gamma = 2$, then it is clear that \eqref{A2} becomes equality. So let us assume $\gamma > 2$ and set
\begin{equation*}
f(x) := (1 - x)^{\gamma} -1 + \gamma x - \frac{\gamma(\gamma - 1)}{2} x^2,
\qquad
x \in [0, 1].
\end{equation*}
Then,
\begin{equation*}
f'(x) = -\gamma (1 - x)^{\gamma -1} + \gamma - \gamma(\gamma - 1) x
=  -\gamma \left[(1 - x)^{\gamma - 1} - 1 + (\gamma-1) x\right].
\end{equation*}
Thus, in view of \eqref{A1} and the fact that $\gamma - 1 > 1$, we get that $f'(x) \leq 0$ for $x \in [0, 1]$, with equality only in the case where $x = 0$. It follows that $f(x)$ is strictly decreasing on $[0, 1]$ and, consequently,
\begin{equation*}
f(x) < f(0) = 0,
\qquad
x \in (0, 1].
\end{equation*}
\hfill $\blacksquare$

\medskip

From Lemma 1 it follows that
\begin{equation}
0 \leq (1 - x)^{\gamma} - (1 - \gamma x) \leq \frac{\gamma(\gamma - 1)}{2} x^2
\qquad
\text{whenever}\ \ x \in [0, 1] \ \ \text{and}\ \ \gamma \geq 2,
\label{A2a}
\end{equation}

\subsection{Asymptotics of the expectation of $T(N)$}

\textbf{Proposition 1.} If for a fixed $a > 0$ we have
\begin{equation}
p \sim \frac{a}{N},
\qquad
N \to \infty,
\label{A0}
\end{equation}
then
\begin{equation}
\mathbb{E}[T(N)] \sim \frac{3a}{2}\log_2 N,
\qquad
N \to \infty,
\label{A3}
\end{equation}
where $\log_2 N$ is the binary logarithm of $N$, namely the logarithm of $N$ to the base $2$ (of course, $\log_2 N = \ln N / \ln 2$).

\smallskip

\textit{Proof}. Let us first consider the case where $N = 2^n$. It was established in \cite{P} that
\begin{equation}
\mathbb{E}\left[T(2^n)\right] = 2^{n+1} - 1 - 2^n\sum_{k=1}^n \frac{q^{2^k} + q^{2^{k-1}}}{2^k}
\qquad
(q = 1 - p).
\label{A4}
\end{equation}
Since
\begin{equation}
\sum_{k=1}^n \frac{q^{2^{k-1}}}{2^k} = \sum_{k=0}^{n-1} \frac{q^{2^k}}{2^{k+1}} 
= \frac{q}{2} + \frac{1}{2} \sum_{k=1}^{n-1} \frac{q^{2^k}}{2^k},
\label{A5}
\end{equation}
formula \eqref{A4} can be written equivalently as
\begin{equation}
\mathbb{E}\left[T(2^n)\right] = 3 \cdot 2^{n-1}\left(1 - \sum_{k=1}^{n-1} \frac{q^{2^k}}{2^k}\right) + 2^{n-1} (1 - q) - q^{2^n} - 1.
\label{A6}
\end{equation}
Now, in view of \eqref{A0} we have that $q = 1 - (a + \varepsilon_N)/N = 1 - (a + \varepsilon_N) 2^{-n}$, where $\varepsilon_N \to 0$ as $N \to \infty$.
Thus, by \eqref{A2a} (with $x = (a + \varepsilon_N) 2^{-n}$ and $\gamma = 2^k$) we get that
\begin{equation}
0 \leq q^{2^k} - 1 + (a + \varepsilon_N) 2^{k-n} \leq
(a + \varepsilon_N)^2 2^{k-1} \left(2^k - 1 \right) 2^{-2n},
\label{A7}
\end{equation}
hence by dividing by $2^k$ and summing from $k=1$ to $k = n-1$ we get
\begin{equation*}
0 \leq \sum_{k=1}^{n-1} \frac{q^{2^k}}{2^k} - \sum_{k=1}^{n-1} \frac{1}{2^k} + (a + \varepsilon_N) (n-1) 2^{-n}
\leq
\frac{(a + \varepsilon_N)^2}{2} 2^{-2n} \sum_{k=1}^{n-1} \left(2^k - 1 \right)
\end{equation*}
or
\begin{equation*}
0 \leq \sum_{k=1}^{n-1} \frac{q^{2^k}}{2^k} - 1 + 2^{1-n} + (a + \varepsilon_N) (n-1) 2^{-n}
\leq
(a + \varepsilon_N)^2 2^{-(2n+1)}  \left(2^n - n - 1 \right).
\end{equation*}
It follows that
\begin{equation}
\sum_{k=1}^{n-1} \frac{q^{2^k}}{2^k} = 1 - (a + \varepsilon_N) n 2^{-n} + O\left(2^{-n}\right),
\qquad
n \to \infty
\label{A10}
\end{equation}
and formula \eqref{A3} follows by substituting \eqref{A10} in \eqref{A6} and using the facts that $q = 1 + O(2^{-n})$ and $n = \log_2 N$.

For the case of a general $N$, i.e. not necessarily a power of $2$, we have, in view of \eqref{I4}, that
\begin{equation}
T\left(2^{\lfloor \nu \rfloor}\right) \leq_{st} T(N) \leq_{st} T\left(2^{\lceil \nu \rceil}\right),
\qquad
\text{where}\ \ 
\nu := \log_2 N,
\label{A11}
\end{equation}
thus
\begin{equation}
\mathbb{E}\left[T\left(2^{\lfloor \nu \rfloor}\right)\right] \leq \mathbb{E}\left[T\left(N\right)\right]
\leq \mathbb{E}\left[T\left(2^{\lceil \nu \rceil}\right)\right].
\label{A11a}
\end{equation}
Now, from our previous analysis we have
\begin{equation}
\mathbb{E}\left[T\left(2^{\lfloor \nu \rfloor}\right)\right] \sim \frac{3a}{2} \lfloor \nu \rfloor,
\qquad
N \to \infty,
\label{A12}
\end{equation}
and
\begin{equation}
\mathbb{E}\left[T\left(2^{\lceil \nu \rceil}\right)\right] \sim \frac{3a}{2} \lceil \nu \rceil,
\qquad
N \to \infty.
\label{A13}
\end{equation}
Therefore, since $\lfloor \nu \rfloor \sim \log_2 N \sim \lceil \nu \rceil$,
formula \eqref{A3} follows by using \eqref{A12} and \eqref{A13} in \eqref{A11a}.
\hfill $\blacksquare$

\medskip

In view of \eqref{I5}, Proposition 1 implies that if $p \ll 1/N \, (p \gg 1/N)$, then $\mathbb{E}[T(N)] \ll \ln N \, (\mathbb{E}[T(N)] \gg \ln N)$.

\medskip


\textbf{Proposition 1A.} If
\begin{equation}
p = \frac{a}{N},
\qquad
a > 0,
\label{A0a}
\end{equation}
then
\begin{equation}
\mathbb{E}[T(N)] = \frac{3a}{2}\log_2 N + O(1),
\qquad
N \to \infty.
\label{A3a}
\end{equation}

\smallskip

\textit{Proof}. Again we first first consider the case where $N = 2^n$. In view of \eqref{A0a} we have that $q = 1 - (a/N) = 1 - a 2^{-n}$.
Thus, by \eqref{A2a} (with $\gamma = 2^k$ and $x = a 2^{-n}$) we get that
\begin{equation}
0 \leq q^{2^k} - 1 + a 2^{k-n} \leq
a^2 2^{k-1} \left(2^k - 1 \right) 2^{-2n},
\label{A7a}
\end{equation}
hence by dividing by $2^k$ and summing from $k=1$ to $k = n-1$ we get
\begin{equation*}
0 \leq \sum_{k=1}^{n-1} \frac{q^{2^k}}{2^k} - \sum_{k=1}^{n-1} \frac{1}{2^k} + a (n-1) 2^{-n}
\leq
\frac{a^2}{2} 2^{-2n} \sum_{k=1}^{n-1} \left(2^k - 1 \right)
\end{equation*}
or
\begin{equation*}
0 \leq \sum_{k=1}^{n-1} \frac{q^{2^k}}{2^k} - 1 + 2^{1-n} + a (n-1) 2^{-n}
\leq
a^2 2^{-(2n+1)}  \left(2^n - n - 1 \right).
\end{equation*}
It follows that
\begin{equation}
\sum_{k=1}^{n-1} \frac{q^{2^k}}{2^k} = 1 - a n 2^{-n} + O\left(2^{-n}\right),
\qquad
n \to \infty
\label{A10a}
\end{equation}
and formula \eqref{A3a} follows by substituting \eqref{A10a} in \eqref{A6} and using the facts that $q = 1 + O(2^{-n})$ and $n = \log_2 N$.

For the case of a general $N$ we have from our previous analysis that
\begin{equation}
\mathbb{E}\left[T\left(2^{\lfloor \nu \rfloor}\right)\right] = \frac{3a}{2} \lfloor \nu \rfloor + O(1),
\qquad
N \to \infty,
\label{A12a}
\end{equation}
and
\begin{equation}
\mathbb{E}\left[T\left(2^{\lceil \nu \rceil}\right)\right] = \frac{3a}{2} \lceil \nu \rceil + O(1),
\qquad
N \to \infty.
\label{A13a}
\end{equation}
Therefore, since $\lfloor \nu \rfloor \leq \log_2 N \leq \lceil \nu \rceil$ and $0 \leq \lceil \nu \rceil - \lfloor \nu \rfloor \leq 1$,
formula \eqref{A3a} follows by using \eqref{A12a} and \eqref{A13a} in \eqref{A11a}.
\hfill $\blacksquare$

\medskip


\textbf{Proposition 2.} If
\begin{equation}
p \sim \frac{a}{N^{\beta}},
\qquad
a > 0, \ \beta > 1,
\label{A14}
\end{equation}
then
\begin{equation}
\mathbb{E}[T(N)] = 1 + O\left(\frac{\ln N}{N^{\beta-1}}\right),
\qquad
N \to \infty.
\label{A15}
\end{equation}
Consequently,
\begin{equation}
\mathbb{E}[T(N) - 1] \to 0,
\qquad
N \to \infty,
\label{A16}
\end{equation}
hence (since $T(N) \geq 1$) the random variable $T(N)$ approaches $1$ in the $L_1$-sense as $N \to \infty$ (actually, from Proposition 5 below it follows that the 
convergence is in the $L_2$-sense).

\smallskip

\textit{Proof}. As in the proof of Proposition 1, let us first consider the case where $N = 2^n$.

In view of \eqref{A14} we have that $q = 1 - (a + \varepsilon_N)/N^{\beta} = 1 - (a + \varepsilon_N) 2^{-\beta n}$, where $\varepsilon_N \to 0$ as $N \to \infty$. 
Thus, by invoking \eqref{A2a} we get that
\begin{equation}
0 \leq q^{2^k} - 1 + (a + \varepsilon_N) 2^{k - \beta n} \leq
(a + \varepsilon_N)^2 2^{k-1} \left(2^k - 1 \right) 2^{-2\beta n},
\label{A17}
\end{equation}
hence by dividing by $2^k$ and summing from $k=1$ to $k = n-1$ we get
\begin{equation}
0 \leq \sum_{k=1}^{n-1} \frac{q^{2^k}}{2^k} - 1 + 2^{1-n} + (a + \varepsilon_N) (n-1) 2^{-\beta n}
\leq
(a + \varepsilon_N)^2 2^{-(2\beta n+1)}  \left(2^n - n - 1 \right).
\label{A17a}
\end{equation}
It follows that (since $\beta > 1$)
\begin{equation}
\sum_{k=1}^{n-1} \frac{q^{2^k}}{2^k} = 1 - 2^{1-n} - (a + \varepsilon_N) (n-1) 2^{-\beta n} + O\left(2^{-(2\beta -1) n}\right),
\qquad
n \to \infty.
\label{A18}
\end{equation}
Let us also notice that by using \eqref{A2a} in \eqref{A14} (with $N = 2^n$) we obtain
\begin{equation}
q^{2^n} = 1 - (a + \varepsilon_N) 2^{-(\beta-1) n} + O\left(2^{-(2\beta -1) n}\right),
\qquad
n \to \infty.
\label{A19}
\end{equation}
Substituting \eqref{A18} and \eqref{A19} in \eqref{A6} and noticing that $q = 1 + O(2^{-\beta n})$ we arrive at the asymptotic formula
\begin{equation}
\mathbb{E}[T(2^n)] = 1 + \frac{3(a + \varepsilon_N)}{2} n 2^{-(\beta-1)n} + O\left(2^{-2(\beta-1)n}\right),
\qquad
n \to \infty,
\label{A20}
\end{equation}
which implies \eqref{A15} in the special case where $N = 2^n$.

In the case of a general $N$ formula \eqref{A15} can be proved by using \eqref{A20} in \eqref{A11a} (as we did in the proof of Proposition 1).
\hfill $\blacksquare$

\medskip

\textbf{Remark 1.} By formula \eqref{A15} and Markov's inequality we get that
\begin{equation}
\mathbb{P}\left\{T(N) - 1 > \frac{1}{\ln N}\right\} = O\left(\frac{(\ln N)^2}{N^{\beta-1}}\right),
\qquad
N \to \infty.
\label{A20a}
\end{equation}
It follows that if $\beta > 2$, so that $(\ln N)^2 /N^{\beta-1}$ is summable, the 1st Borel-Cantelli Lemma and the fact that $T(N) \geq 1$ imply that,
under \eqref{A14}, $T(N) \to 1$ a.s. as $N \to \infty$.

\medskip

\textbf{Proposition 3.} If
\begin{equation}
p \sim \frac{a}{N^{\beta}},
\qquad
a > 0, \ 0 < \beta < 1,
\label{A21}
\end{equation}
then for every $\delta > 0$ there is an $N_0 = N_0(\delta)$ such that
\begin{equation}
\mathbb{E}[T(N)] \leq 2^{1 - \beta} \left(\frac{3a}{2} + \delta\right) N^{1-\beta} \log_2 N
\qquad \text{for all}\ \ 
N \geq N_0.
\label{A22}
\end{equation}

\smallskip

\textit{Proof}. Again we first considet the case $N = 2^n$.

Formula \eqref{A17a} still holds under \eqref{A21}, namely for $\beta \in (0, 1)$ and, thus, by the first inequality of \eqref{A17a} we get
\begin{equation}
1 - \sum_{k=1}^{n-1} \frac{q^{2^k}}{2^k} \leq 2^{1-n} + (a + \varepsilon_N) (n-1) 2^{-\beta n}.
\label{A23}
\end{equation}
Now, in view of \eqref{A23} and \eqref{A21}, formula \eqref{A6} yields
\begin{equation*}
\mathbb{E}\left[T(2^n)\right] \leq \frac{3(a + \varepsilon_N)}{2} n 2^{(1-\beta) n} - (a + \varepsilon_N) 2^{(1-\beta) n} + 2 - q^{2^n}.
\end{equation*}
It follows that for any given $\delta > 0$ there is an $n_0 = n_0(\delta)$ such that
\begin{equation}
\mathbb{E}\left[T(2^n)\right] \leq \left(\frac{3a}{2} + \delta\right) n 2^{(1-\beta) n},
\label{A24}
\end{equation}
and this implies \eqref{A22} in the special case where $N = 2^n$.

Finally, the case of a general $N$ follows from \eqref{A24} together with the second inequality of formula \eqref{A11a}, namely the inequality
\begin{equation*}
\mathbb{E}\left[T\left(N\right)\right]
\leq \mathbb{E}\left[T\left(2^{\lceil \nu \rceil}\right)\right],
\qquad \text{where} \ \ 
\nu = \log_2 N.
\end{equation*}
\hfill $\blacksquare$

\medskip

Notice that Proposition 3 tells us that, even under \eqref{A21}, our binary search scheme still behaves much better than individual testing, at least for sufficiently
large $N$.

\subsection{Asymptotics of the variance of $T(N)$}

In \cite{P} we had derived an explicit formula for the variance $\mathbb{V}[T(N)]$ of $T(N)$ in the case where $N = 2^n$:
\begin{equation}
\mathbb{V}[T(2^n)] = V_1 + V_2,
\label{A26a}
\end{equation}
where
\begin{equation}
V_1 := 2^n \sum_{k=1}^n \left(2 q^{2^k} + q^{2^{k-1}}\right) \left(2 -  \sum_{j=1}^k  \frac{q^{2^j} + q^{2^{j-1}}}{2^j}\right)
\label{A26b}
\end{equation}
and
\begin{equation}
 V_2 := 2^n \sum_{k=1}^n \frac{q^{2^{k+1}} + q^{3 \cdot 2^{k-1}} - 5 q^{2^k} - 3 q^{2^{k-1}}}{2^k}
\label{A26c}
\end{equation}
(in the trivial case $n=0$ all the above sums are empty, i.e. $0$).

By using the formula
\begin{equation*}
2 = \frac{1}{2^{k-1}} + \sum_{j=1}^k  \frac{2}{2^j}
\end{equation*}
in the second parentheses of the right-hand side of \eqref{A26b}, we obtain
\begin{align}
V_1 &= 2^n \sum_{k=1}^n \left(2 q^{2^k} + q^{2^{k-1}}\right) \left(\frac{1}{2^{k-1}} + \sum_{j=1}^k  \frac{2 - q^{2^j} - q^{2^{j-1}}}{2^j}\right)
\nonumber
\\
&= V_3 + 2^n \sum_{k=1}^n \frac{4 q^{2^k} + 2q^{2^{k-1}}}{2^k},
\label{A27}
\end{align}
where
\begin{equation}
 V_3 := 2^n \sum_{k=1}^n \left(2 q^{2^k} + q^{2^{k-1}}\right) \sum_{j=1}^k  \frac{2 - q^{2^j} - q^{2^{j-1}}}{2^j}.
\label{A27a}
\end{equation}
Also, in view of \eqref{A26c},
\begin{equation}
V_4 := V_2 + 2^n \sum_{k=1}^n \frac{4 q^{2^k} + 2q^{2^{k-1}}}{2^k} =  2^n \sum_{k=1}^n \frac{q^{2^{k+1}} + q^{3 \cdot 2^{k-1}} - q^{2^k} - q^{2^{k-1}}}{2^k},
\label{A27b}
\end{equation}
so that
\begin{equation}
\mathbb{V}[T(2^n)] = V_3 + V_4.
\label{A26d}
\end{equation}


\medskip

\textbf{Proposition 4.} Suppose $p$ satisfies \eqref{A0a}, namely $p = a /N$. Then
\begin{equation}
\mathbb{V}[T(N)] = \frac{9a}{4} (\log_2 N)^2 + O(\ln N),
\qquad
N \to \infty.
\label{A28}
\end{equation}

\smallskip

\textit{Proof}. Let us first consider the case where $N = 2^n$.

By applying \eqref{A7} we get
\begin{equation}
-\frac{a^2}{4} 2^{-n} \left(2^{j+1} + 2^{j-1} - 3\right) \leq 2^n \frac{2 - q ^{2^j} - q ^{2^{j-1}}}{2^j} - \frac{3a}{2} \leq 0.
\label{A29}
\end{equation}
Taking the double summation $\sum_{k=1}^n \left(2 q^{2^k} + q^{2^{k-1}}\right) \sum_{j=1}^k$ to all terms of \eqref{A29} yields, in view of \eqref{A27a}, 
\begin{align}
-\frac{a^2}{4} 2^{-n} \sum_{k=1}^n &\left(2 q^{2^k} + q^{2^{k-1}}\right) \sum_{j=1}^k \left(2^{j+1} + 2^{j-1} - 3\right)
\nonumber
\\
&\leq V_3 - \frac{3a}{2} \sum_{k=1}^n \left(2 q^{2^k} + q^{2^{k-1}}\right) k
\leq 0.
\label{A30}
\end{align}
Now, for $1 \leq k \leq n$ we have
\begin{equation*}
2^{-n} \sum_{j=1}^k \left(2^{j+1} + 2^{j-1} - 3\right) \leq 2^{-n} \sum_{j=1}^n 2^{j+2} \leq 8,
\end{equation*}
hence the left-hand side of \eqref{A30} is (at most) $O(n)$ as $n \to \infty$. Therefore, \eqref{A30} implies
\begin{equation}
V_3 = \frac{3a}{2} \sum_{k=1}^n \left(2 q^{2^k} + q^{2^{k-1}}\right) k + O(n),
\qquad
n \to \infty.
\label{A32}
\end{equation}
Next, we estimate the sum in the right-hand side of \eqref{A32}. Again by \eqref{A7} we get
\begin{equation}
0 \leq 
\left(2 q^{2^k} + q^{2^{k-1}}\right) k - 3k + 5a 2^{k-1-n} k \leq
a^2 2^{k-2-2n} \left(9 \cdot 2^{k-1} - 5\right) k.
\label{A33}
\end{equation}
Summing all terms of \eqref{A33} from $k=1$ to $k=n$ and recalling that
\begin{equation}
\sum_{k=1}^n \rho^{k-1} k = \frac{(\rho-1)\rho^n n - \rho^n + 1 }{(\rho-1)^2},
\qquad
\rho \ne 1,
\label{A34}
\end{equation}
we get that
\begin{equation*}
\sum_{k=1}^n \left(2 q^{2^k} + q^{2^{k-1}}\right) k = \sum_{k=1}^n 3k + O(n) = \frac{3}{2} n^2 + O(n),
\qquad
n \to \infty.
\end{equation*}
Therefore, \eqref{A32} yields
\begin{equation}
V_3 = \frac{9a}{4} n^2 + O(n),
\qquad
n \to \infty.
\label{A36}
\end{equation}
Finally, we need to estimate the quantity $V_4$ of \eqref{A27b}. In the same way we derived formula \eqref{A10} we can obtain the estimates
\begin{equation}
\sum_{k=1}^n \frac{q^{2^{k+1}}}{2^k} = 1 - 2a n 2^{-n} + O\left(2^{-n}\right),
\qquad
n \to \infty,
\label{A37a}
\end{equation}
\begin{equation}
\sum_{k=1}^n \frac{q^{3\cdot 2^{k-1}}}{2^k} = 1 - \frac{3a}{2} n 2^{-n} + O\left(2^{-n}\right),
\qquad
n \to \infty,
\label{A37b}
\end{equation}
\begin{equation}
\sum_{k=1}^n \frac{q^{2^k}}{2^k} = 1 - a n 2^{-n} + O\left(2^{-n}\right),
\qquad
n \to \infty,
\label{A37c}
\end{equation}
and
\begin{equation}
\sum_{k=1}^n \frac{q^{2^{k-1}}}{2^k} = 1 - \frac{a}{2} n 2^{-n} + O\left(2^{-n}\right),
\qquad
n \to \infty.
\label{A37d}
\end{equation}
Using \eqref{A37a}, \eqref{A37b}, \eqref{A37c}, and \eqref{A37d} in \eqref{A27b} we obtain
\begin{equation}
V_4 = -2an + O(1),
\qquad
n \to \infty.
\label{A38}
\end{equation}
Hence, in view of \eqref{A38} and \eqref{A36}, formula \eqref{A26d} yields
\begin{equation}
\mathbb{V}[T(2^n)] = \frac{9a}{4} n^2 + O(n),
\qquad
n \to \infty.
\label{A39}
\end{equation}
which verifies \eqref{A28} in the case where $N = 2^n$. Furthermore, from \eqref{A39} and \eqref{A3} we have that
\begin{equation}
\mathbb{E}\left[T(2^n)^2\right] = \mathbb{V}[T(2^n)] + \mathbb{E}[T(2^n)]^2 = \frac{9a(a+1)}{4} n^2 + O(n),
\qquad
n \to \infty.
\label{A40}
\end{equation}

For the case of a general $N$, starting from \eqref{A40} we can use the same approach we used in the proof of Proposition 1 to show that
\begin{equation}
\mathbb{E}\left[T(N)^2\right] = \frac{9a(a+1)}{4} (\log_2 N)^2 + O(\ln N),
\qquad
N \to \infty.
\label{A41}
\end{equation}
Then, \eqref{A28} follows from \eqref{A41} and \eqref{A3}.
\hfill $\blacksquare$

\medskip

By adapting the proof of Proposition 4 (in the spirit of the proof of Proposition 1) we can show that if $p \sim a /N$, then
\begin{equation*}
\mathbb{V}[T(N)] \sim \frac{9a}{4} (\log_2 N)^2,
\qquad
N \to \infty.
\end{equation*}

\medskip

\textbf{Proposition 5.} Suppose $p = a /N^{\beta}$ for some fixed $a > 0$ and $\beta > 1$.
Then
\begin{equation}
\mathbb{V}[T(N)] = O\left(\frac{\ln^2 N}{N^{\beta-1}}\right),
\qquad
N \to \infty.
\label{B28}
\end{equation}
Consequently, in comparison with Proposition 2, the random variable $T(N)$ approaches $1$ in the $L_2$-sense as $N \to \infty$.

\smallskip

\textit{Proof}. As usual, we first consider the case where $N = 2^n$.

In view of our assumption for $p$ we have that $q = 1 - (a/N^{\beta}) = 1 - a 2^{-\beta n}$. Thus, by invoking \eqref{A2a} we get that
\begin{equation}
-\frac{a^2}{4} 2^{-(2\beta-1)n} \left(2^{j+1} + 2^{j-1} - 3\right) \leq 2^n \frac{2 - q ^{2^j} - q ^{2^{j-1}}}{2^j} - \frac{3a}{2} 2^{-(\beta-1)n} \leq 0.
\label{B29}
\end{equation}
Taking the double summation $\sum_{k=1}^n \left(2 q^{2^k} + q^{2^{k-1}}\right) \sum_{j=1}^k$ to all terms of \eqref{B29} yields, in view of \eqref{A27a}, 
\begin{align}
-\frac{a^2}{4} 2^{-(2\beta-1)n} \sum_{k=1}^n &\left(2 q^{2^k} + q^{2^{k-1}}\right) \sum_{j=1}^k \left(2^{j+1} + 2^{j-1} - 3\right)
\nonumber
\\
&\leq V_3 - \frac{3a}{2} 2^{-(\beta-1)n} \sum_{k=1}^n \left(2 q^{2^k} + q^{2^{k-1}}\right) k
\leq 0.
\label{B30}
\end{align}
As we have already noticed in the proof of Proposition 4, for $1 \leq k \leq n$ we have
\begin{equation*}
2^{-n} \sum_{j=1}^k \left(2^{j+1} + 2^{j-1} - 3\right) \leq 2^{-n} \sum_{j=1}^n 2^{j+2} \leq 8,
\end{equation*}
hence the left-hand side of \eqref{B30} is (at most) $O\left(2^{-2(\beta-1)n} n\right)$ as $n \to \infty$. Therefore, \eqref{B30} implies
\begin{equation}
V_3 = \frac{3a}{2} 2^{-(\beta-1)n} \sum_{k=1}^n \left(2 q^{2^k} + q^{2^{k-1}}\right) k + O\left(2^{-2(\beta-1)n} n\right),
\qquad
n \to \infty.
\label{B32}
\end{equation}
Next, we estimate the sum in the right-hand side of \eqref{B32}. Again by \eqref{A2a} we get
\begin{equation}
0 \leq 
\left(2 q^{2^k} + q^{2^{k-1}}\right) k - 3k + \frac{5a}{2} 2^{k - \beta n} k \leq
\frac{a^2}{4} 2^{k - 2\beta n} \left(9 \cdot 2^{k-1} - 5\right) k.
\label{B33}
\end{equation}
Summing all terms of \eqref{B33} from $k=1$ to $k=n$ and recalling \eqref{A34} we get that, as $n \to \infty$,
\begin{equation*}
\sum_{k=1}^n \left(2 q^{2^k} + q^{2^{k-1}}\right) k = \sum_{k=1}^n 3k + O\left(2^{-(\beta-1)n} n\right) = \frac{3}{2} n^2 + O\left(2^{-(\beta-1)n} n\right).
\end{equation*}
Therefore, \eqref{B32} yields
\begin{equation}
V_3 = \frac{9a}{4} 2^{-(\beta-1)n} n^2 + O\left(2^{-(\beta-1)n} n\right),
\qquad
n \to \infty.
\label{B36}
\end{equation}

Finally, we need to estimate the quantity $V_4$ of \eqref{A27b}. In the same way we derived formula \eqref{A10} we can obtain the estimates
\begin{equation}
\sum_{k=1}^n \frac{q^{2^{k+1}}}{2^k} = 1 - 2^{-n} - 2a n 2^{-\beta n} + O\left(2^{-(2\beta - 1) n}\right),
\qquad
n \to \infty,
\label{B37a}
\end{equation}
\begin{equation}
\sum_{k=1}^n \frac{q^{3\cdot 2^{k-1}}}{2^k} = 1 - 2^{-n} - \frac{3a}{2} n 2^{-\beta n} + O\left(2^{-(2\beta - 1) n}\right),
\qquad
n \to \infty,
\label{B37b}
\end{equation}
\begin{equation}
\sum_{k=1}^n \frac{q^{2^k}}{2^k} = 1 - 2^{-n} - a n 2^{-\beta n} + O\left(2^{-(2\beta - 1) n}\right),
\qquad
n \to \infty,
\label{B37c}
\end{equation}
and
\begin{equation}
\sum_{k=1}^n \frac{q^{2^{k-1}}}{2^k} = 1 - 2^{-n} - \frac{a}{2} n 2^{-\beta n} + O\left(2^{-(2\beta - 1) n}\right),
\qquad
n \to \infty.
\label{B37d}
\end{equation}
Using \eqref{B37a}, \eqref{B37b}, \eqref{B37c}, and \eqref{B37d} in \eqref{A27b} we obtain
\begin{equation}
V_4 = -2a 2^{-(\beta - 1)n} n + O\left(2^{-2(\beta - 1)n}\right),
\qquad
n \to \infty.
\label{B38}
\end{equation}
Hence, in view of \eqref{B38} and \eqref{B36}, formula \eqref{A26d} yields
\begin{equation}
\mathbb{V}[T(2^n)] =  \frac{9a}{4} 2^{-(\beta-1)n} n^2 + O\left(2^{-(\beta-1)n} n\right),
\qquad
n \to \infty.
\label{B39}
\end{equation}
which verifies \eqref{B28} in the case where $N = 2^n$. Furthermore, from \eqref{B39} and \eqref{A20} we have, as $n \to \infty$,
\begin{align}
\mathbb{E}\left[T(2^n)^2\right] &= \mathbb{V}[T(2^n)] + \mathbb{E}[T(2^n)]^2
\nonumber
\\
&= 1 + \frac{9a}{4} 2^{-(\beta-1)n} n^2 + 3a 2^{-(\beta-1)n} n + O\left(2^{-2(\beta-1)n} n^2\right).
\label{B40}
\end{align}

For the case of a general $N$, starting from \eqref{B40} we can use the same approach we used in the proof of Proposition 1 to show that
\begin{equation}
\mathbb{E}\left[T(N)^2\right] = 1 + O\left(\frac{\ln^2 N}{N^{\beta - 1}}\right),
\qquad
N \to \infty.
\label{B41}
\end{equation}
Then, \eqref{B28} follows from \eqref{B41} and \eqref{A15}.
\hfill $\blacksquare$

\medskip

The statement of Proposition 5 remains true (with essentially the same proof) under the more general assumption that $p \sim a /N^{\beta}$, $N \to \infty$,
for some fixed $a > 0$ and $\beta > 1$.

\medskip

Finally, let us propose an open question.

\medskip

\textbf{Open Question.} Suppose $p = a/N$. Determine the limiting distribution(s) of
\begin{equation*}
\frac{T(N)}{\ln N}
\qquad
\end{equation*}
as $N \to \infty$. Notice that $T(N)/\ln N > 0$, hence the limiting distribution cannot, in particular, be normal (in contrast to the case of a fixed $p$ where the
limiting distribution of $T(N)$, appropriately normalized, is normal \cite{P}).


\end{document}